Research

# On Godunov-type finite volume methods for seismic wave propagation

Juan Barrios[1] · Pedro S. Peixoto[1] · Felipe A. G. Silva[1]

**Abstract**
The computational complexity of simulating seismic waves demands continual exploration of more efficient numerical methods. While Finite Volume methods are widely acclaimed for tackling general nonlinear hyperbolic (wave) problems, their application in realistic seismic wave simulation remains uncommon, with rare investigations in the literature. Furthermore, seismic wavefields are influenced by sharp subsurface interfaces frequently encountered in realistic models, which could, in principle, be adequately solved with Finite Volume methods. In this study, we delved into two Finite Volume (FV) methods to assess their efficacy and competitiveness in seismic wave simulations, compared to traditional Finite Difference schemes. We investigated Gudunov-type FV methods: an upwind method called wave propagation algorithm (WPA), and a Central-Upwind type method (CUp). Our numerical analysis uncovered that these finite volume methods could provide less dispersion (albeit increased dissipation) compared to finite differences for seismic problems characterized by velocity profiles with abrupt transitions in the velocity. However, when applied to more realistic seismic models, finite volume methods yielded unfavorable outcomes compared to finite difference methods, the latter offering lower computational costs and higher accuracy. This highlights that despite the potential advantages of finite volume methods, such as their conservative nature and aptitude for accurately capturing shock waves in specific contexts, our results indicate that they are only advantageous for seismic simulations when unrealistic abrupt transitions are present in the velocity models.

## 1 Introduction

Seismic wave propagation, a cornerstone of computational seismology, involves numerically solving elastic wave propagation problems. This process, often referred to as the direct problem in seismic imaging, lies at the core of various applications within computational seismology. These applications span a broad spectrum, including geophysical exploration, regional wave propagation analysis, global or planetary seismology studies, investigations into ground motion and dynamic rupture phenomena, seismic topography (utilizing techniques like full waveform inversion), among others [1].

Different mathematical models have been developed to simulate wave propagation in the complex geological structure of the earth, using various numerical techniques, such as finite differences, finite elements, finite volumes, spectral or boundary element methods, discontinuous Galerkin method, and hybrid methods, which have been constantly improved to obtain an accurate and reliable solution of the wave equation for various media [1, 2]. Still, due to its low computational cost and simplicity, the finite difference method remains the most widely used method. However, although the finite difference method is currently found in large applications in the analysis of seismic wave propagation, in complex earth structures, including those with large velocity contrasts, strong heterogeneity, relief, or topographic attenuation, the methodology can fail [3]. This is mainly due to the assumption often used in formulating finite difference methods that both the wave propagation medium and the solution function of the problem are smooth functions. For problems with discontinuities, the accuracy of finite difference methods deteriorates and may fail to represent propagation through regions of abrupt transitions.

✉ Juan Barrios, jbarrioscamargo@usp.br; Pedro S. Peixoto, ppeixoto@usp.br; Felipe A. G. Silva, felipe.augusto.guedes@usp.br | [1]Instituto de Matemática e Estatística, Universidade de São Paulo, São Paulo, SP, Brazil.

The finite volume method offers a great advantage by generally not requiring smoothness of the propagation medium or the solution and is a natural choice for modeling in heterogeneous media. In this type of method, each mesh cell can be given different material properties. The propagation of information between cells is done by analyzing the physics imposed on the cell boundaries, and is suitable and popular for hyperbolic conservation laws and equilibrium law systems; See, for example, [4–10].

In this numerical study, we investigated two Godunov-type finite volume methods, based on the **REA** algorithm (Reconstruct—Evolve—Average), to explore their potential in modeling seismic waves across heterogeneous media.

The first method, known as the Wave Propagation Algorithm (WPA), is an Upwind type, originally developed by [11]. This method leverages the resulting waves from each Riemann solution in conjunction with limiter flux functions to achieve sharp resolution of discontinuities without inducing numerical oscillations. It exhibits second-order accuracy where the solution remains smooth and has demonstrated efficacy in numerically solving various wave propagation problems, for example, advection equation, linear acoustics equations (in homogeneous and heterogeneous media), Burger's equation, Euler equations for gas dynamics, shallow water equations, and even the elasticity equations, as well as examples of applications in tsunami simulation.

The second method is a Central-Upwind (CUp) type method pioneered by [12], building upon the one-dimensional fully and semi-discrete central-upwind schemes introduced in earlier work [13]. This method aims to reduce numerical dissipation, providing a more accurate solution of the evolved quantities on the original grid. The scheme's simplicity and generality make it an attractive alternative for solving wave propagation problems.

Its applications can be seen, for example: in [12] where the scheme is applied to solve Euler equations; in [8] a central-upwind method for shallow water equations is developed and, most recently, in [14] an improvement on the numerical dissipation of Central-Upwind schemes with implementation on a wide variety of complex numerical examples.

In the existing literature, very few different finite volume methods have been applied to seismic wave propagation or full waveform inversion problems. Some notable articles address the problem for unstructured grids. For instance, in [15], they present the development of arbitrary high-order finite volume methods on unstructured meshes for the seismic wave propagation problem in 2D and 3D settings. [16] proposed a comprehensive reanalysis of the finite-volume approach based on unstructured triangular meshes. More recently, [17] presented a cell-centered finite volume scheme for the diffusive-viscous wave equation on general polygonal meshes. Additionally, [18] developed a finite volume method in the frequency domain for 2D problems. A certain common knowledge that Finite Volume methods are expected to be more expensive than Finite Differences for regular cartesian grids is sometimes heard, but rarely documented.

Additionally, while numerical methods for solving elasticity equations, often in their first-order form as coupled hyperbolic systems, have been documented in the literature using techniques such as WPA [19, 20] or modified Central-Upwind methods [21], there is a notable absence of literature that evaluates their performance when applied to seismic problems. This research addresses this important gap in the literature by systematically comparing finite difference and finite volume methods in the context of seismic wave propagation. Our analysis highlights the advantages and limitations of each method, offering novel insights into their performance under varying conditions.

The paper is structured as follows: In Sect. 2, we introduce the equations governing the problem under study. Section 3 provides a detailed exposition of the finite volume methods implemented. Finally, in Sect. 4, we present the numerical results, including their application to synthetic and realistic seismic problems.

## 2 Modeling equations

In this section, we present the equations governing our study problem (for further details see, e.g. [22] or [10]).

### 2.1 2-D elastic wave equations

To model compression waves (P-Waves) in 2-D one can use the following system

$$\begin{aligned}\varepsilon_t - u_x - v_y &= 0, \\ (\rho(x,y)u)_t - \sigma_x(K(x,y);\varepsilon) &= 0, \\ (\rho(x,y)u)_t - \sigma_y(K(x,y);\varepsilon) &= 0,\end{aligned} \tag{1}$$

where, for brevity, the partial derivatives are denoted by subscripts. Here $u$ and $v$ are the velocities in the $x$ and $y$ directions, respectively, $\varepsilon$ is the strain, $\sigma$ is the stress, and $\rho$ is the density.

The system (1) can be written in the form of a conservation law

$$\mathbf{q}_t + f(\mathbf{q}, x, y)_x + g(\mathbf{q}, x, y)_y = 0, \tag{2}$$

subject to the initial condition, $\mathbf{q}(x,0) = \mathbf{q}_0(x)$. Here

$$\mathbf{q} = \begin{pmatrix} \varepsilon \\ m_u \\ m_v \end{pmatrix}, \quad f(\mathbf{q}) = \begin{pmatrix} m_u/\rho \\ -\sigma \\ 0 \end{pmatrix}, \quad g(\mathbf{q}) = \begin{pmatrix} m_v/\rho \\ 0 \\ -\sigma \end{pmatrix}, \tag{3}$$

where $m_u = \rho u$ and $m_v = \rho v$ denote the corresponding momenta.

For small strains, we can assume a linear stress–strain relationship of the form

$$\sigma(\varepsilon, x, y) = K(x,y)\varepsilon,$$

where $K(x, y)$ is the volumetric modulus of compressibility. In this case, the linear hyperbolic system $\mathbf{q}_t + (\mathbf{A}(x,y)\mathbf{q})_x + (\mathbf{B}(x,y)\mathbf{q})_y = 0$ has coefficient matrices

$$A(x,y) = \begin{bmatrix} 0 & -1/\rho(x,y) & 0 \\ -K(x,y) & 0 & 0 \\ 0 & 0 & 0 \end{bmatrix}, B(x,y) = \begin{bmatrix} 0 & 0 & -1/\rho(x,y) \\ 0 & 0 & 0 \\ -K(x,y) & 0 & 0 \end{bmatrix}. \tag{4}$$

This system is a simplification such that only P-waves are modeled, not S-waves [19].

## 3 Numerical methods

In one spatial dimension, the Finite Volume Method (FVM) involves partitioning the domain into discrete intervals known as finite volumes or control volumes and tracking an approximation of the flow integral $\mathbf{q}$ over each of these volumes. At each time step, these values are updated by approximating the flow through the boundaries of the intervals.

The mesh contains control volumes of type $\mathscr{C}_i = (x_{i-1/2}, x_{i+1/2})$, of fixed width $\Delta x = x_{i+1/2} - x_{i-1/2}$, and a time step $\Delta t = t^{n+1} - t^n$, as shown in Fig. 1. The value $\mathbf{Q}_i^n$ will approximate the average value over the $i$-th interval in time $t_n$:

$$\mathbf{Q}_i^n \approx \frac{1}{\Delta x} \int_{x_{i-1/2}}^{x_{i+1/2}} \mathbf{q}(x, t_n) dx \equiv \frac{1}{\Delta x} \int_{\mathscr{C}_i} \mathbf{q}(x, t_n) dx, \tag{5}$$

and

$$\mathbf{F}_{i-1/2}^n = \mathscr{F}(\mathbf{Q}_{i-1}^n, \mathbf{Q}_i^n), \qquad \mathbf{F}_{i+1/2}^n = \mathscr{F}(\mathbf{Q}_i^n, \mathbf{Q}_{i+1}^n), \tag{6}$$

where $\mathscr{F}$ is a numerical flow function on the respective interfaces.

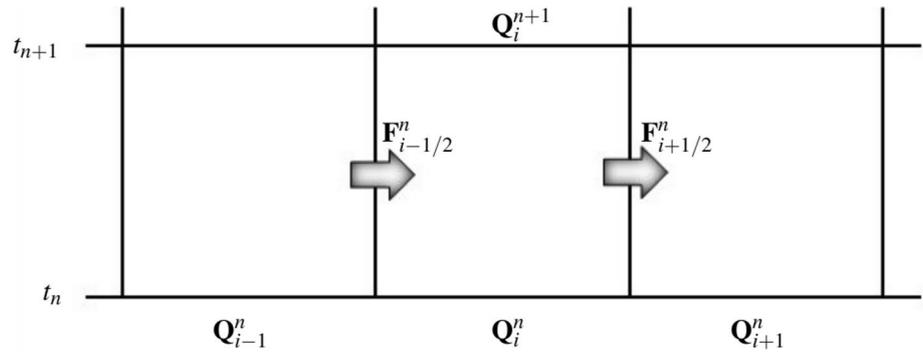

**Fig. 1** Illustration of the variables in a finite volume method to update the value $\mathbf{Q}_i^n$ by the flux at the interfaces of the control volume

## 3.1 High-resolution wave propagation algorithms (WPA)

In this section, we introduce the Wave Propagation Algorithm, an upwind Godunov-type FVM developed by [10]. This algorithm is formulated as a wave propagation scheme for systems of hyperbolic conservation laws and represents a second-order method based on the Lax-Wendroff approach. In fact, it reduces to the classic Lax-Wendroff method in the absence of flux limiters. By incorporating flux limiters (referred to as wave limiters within the context of the Wave Propagation Algorithm), the method achieves high resolution while mitigating non-physical oscillations near discontinuities or steep solution gradients. For a more comprehensive understanding of the Wave Propagation Algorithm, readers are referred to [10].

This method requires solving a Riemann problem at cell edges to update the cell average 5 in each time step. For an autonomous system, the Riemann problem at $x_{i-1/2}$ is a particular case of a Cauchy problem and consists of the hyperbolic equation with constant coefficients,

$$\mathbf{q}_t + \mathbf{A}\mathbf{q}_x = 0, \tag{7}$$

together with a special initial condition that is piecewise constant,

$$\mathbf{q}(x,0) = \begin{cases} \mathbf{Q}_{i-1}, & \text{if } x < x_{i-1/2}, \\ \mathbf{Q}_i, & \text{if } x > x_{i-1/2}. \end{cases} \tag{8}$$

The Riemann solution for a system of $m$ equations typically consists of $m$ waves that we denote by $\mathscr{W}_{i-1/2}^p$, for $p = 1, 2, \cdots, m$, propagating with speeds $s_{i-1/2}^p$. The problem can be easily solved in terms of the eigenvectors $\mathbf{r}_{i-1/2}^p$. The standard approach is to decompose the jump in $\mathbf{Q}$ as a linear combination of the eigenvectors to define waves $\mathscr{W}_{i-1/2}^p$

$$\mathbf{Q}_i - \mathbf{Q}_{i-1} = \sum_{p=1}^m \alpha_{i-1/2}^p \mathbf{r}_{i-1/2}^p \equiv \sum_{p=1}^m \mathscr{W}_{i-1/2}^p. \tag{9}$$

The coefficients $\alpha_{i-1/2}^p$ are given by

$$\alpha_{i1/2}^p = \mathbf{R}_{i-1/2}^{-1}(\mathbf{Q}_i - \mathbf{Q}_{i-1}), \tag{10}$$

where $\mathbf{R}_{i-1/2}$ is the matrix or right eigenvectors.

The **WPA** method takes the form

$$\mathbf{Q}_i^{n+1} = \mathbf{Q}_i - \frac{\Delta t}{\Delta x}(\mathscr{A}^+ \Delta \mathbf{Q}_{i-1/2} + \mathscr{A}^- \Delta \mathbf{Q}_{i+1/2}) - \frac{\Delta t}{\Delta x}(\tilde{F}_{i+1/2} - \tilde{F}_{i-1/2}), \tag{11}$$

where the terms

$$\begin{aligned}
\mathscr{A}^- \Delta \mathbf{Q}_{i-1/2} &= \sum_{p=1}^m (s^p)^- \mathscr{W}_{i-1/2}^p, \\
\mathscr{A}^+ \Delta \mathbf{Q}_{i-1/2} &= \sum_{p=1}^m (s^p)^+ \mathscr{W}_{i-1/2}^p,
\end{aligned} \tag{12}$$

are called fluctuations and the flux correction term is defined by

$$\tilde{F}_{i-1/2} = \frac{1}{2} \sum_{p=1}^{m} |s^p|\left(1 - \frac{\Delta t}{\Delta x}|s^p|\right)\tilde{\mathcal{W}}^p_{i-1/2}, \tag{13}$$

are based on the waves $\mathcal{W}^p_{i-1/2}$ and speeds $s^p_{i-1/2}$ resulting when solving the Riemann problem for any two states $\mathbf{Q}_{i-1}$ and $\mathbf{Q}_i$ and

$$\tilde{\mathcal{W}}^p_{i-1/2} = \psi(\theta^p_{i-1/2})\mathcal{W}^p_{i-1/2}, \tag{14}$$

is a limited version of the original wave, where

$$\theta^p_{i-1/2} = \frac{\mathcal{W}^p_{I-1/2} \cdot \mathcal{W}^p_{i-1/2}}{\mathcal{W}^p_{i-1/2} \cdot \mathcal{W}^p_{i-1/2}} \qquad \text{with } I = \begin{cases} i-1 & \text{if } s^p > 0, \\ i+1 & \text{if } s^p < 0, \end{cases} \tag{15}$$

and $\psi$ is a wave limiter. In this paper, WPA is implemented using the SuperBee (SB) limiter in the flux limiter version

$$\psi(\theta) = \max(0, \min(1, 2\theta), \min(2, \theta)). \tag{16}$$

As the objective of the study is to verify the behavior of this type of method when applied to the elasticity equations in heterogeneous media, we used the f-wave formulation, proposed in [23], which is often convenient for equations with spatially varying flux functions. In this case, the fluctuations are defined as

$$\begin{aligned}\mathcal{A}^- \Delta \mathbf{Q}_{i-1/2} &= \sum_{p:s^p_{i-1/2}<0} \mathcal{Z}^p_{i-1/2}, \\ \mathcal{A}^+ \Delta \mathbf{Q}_{i-1/2} &= \sum_{p:s^p_{i-1/2}>0} \mathcal{Z}^p_{i-1/2},\end{aligned} \tag{17}$$

where we sum only over the $p$ for which $s^p_{i-1/2}$ is negative or positive and replace the correction flow by

$$\tilde{\mathbf{F}}_{i-1/2} = \frac{1}{2} \sum_{p=1}^{M_w} \text{sgn}(s^p_{i-1/2})\left(1 - \frac{\Delta t}{\Delta x}|s^p_{i-1/2}|\right)\tilde{\mathcal{Z}}^p_{i-1/2}, \tag{18}$$

where $\tilde{\mathcal{Z}}^p$ is a limited version of the wave $\mathcal{Z}^p$, that are called **f-waves**, obtained in the same way that $\tilde{\mathcal{W}}^p$ would be obtained from $\mathcal{W}^p$.

In the two-dimensional case, the value $\mathbf{Q}^n_{ij}$ represents an average of cells over the grid cell $(i,j)$ at time $t_n$.

$$\mathbf{Q}_{i,j} \approx \frac{1}{\Delta x \Delta y}\int_{y_{j-1/2}}^{y_{j+1/2}}\int_{x_{i-1/2}}^{x_{i+1/2}} \mathbf{q}(x,y,t_n)dxdy, \tag{19}$$

where $\Delta x = x_{i+1/2,j} - x_{i-1/2,j}$, $\Delta y = y_{i,j+1/2} - y_{i,j-1/2}$.

In the special case of a rectangular grid of the form $\mathcal{C}_{ij} = [x_{i-1/2}, x_{i+1/2}] \times [y_{j-1/2}, y_{j+1/2}]$, as shown in Fig. 2, where $\Delta x = x_{i+1/2} - x_{i-1/2}$ and $\Delta y = y_{j+1/2} - y_{j-1/2}$, we have two different ways to extend the wave propagation algorithm to two dimensions, the first leads us to the **Fully Discrete Wave Propagation Algorithm (FWPA)** given by

$$\begin{aligned}\mathbf{Q}^{n+1}_{ij} =\ & \mathbf{Q}_{ij} - \frac{\Delta t}{\Delta x}(\mathcal{A}^+ \Delta \mathbf{Q}_{i-1/2,j} + \mathcal{A}^- \Delta \mathbf{Q}_{i+1/2,j}) \\ & - \frac{\Delta t}{\Delta x}(\mathcal{B}^+ \Delta \mathbf{Q}_{i,j-1/2} + \mathcal{B}^- \Delta \mathbf{Q}_{i,j+1/2}) \\ & - \frac{\Delta t}{\Delta x}(\tilde{\mathbf{F}}_{i+1/2,j} - \tilde{\mathbf{F}}_{i-1/2,j}) - \frac{\Delta t}{\Delta y}(\tilde{\mathbf{G}}_{i,j+1/2} - \tilde{\mathbf{G}}_{i,j-1/2}),\end{aligned} \tag{20}$$

where $\tilde{\mathbf{F}}$, and $\tilde{\mathbf{G}}$ are calculated via an algorithm that uses information from neighboring cells as described in [10].

The second way leads us to the **Dimensional Splitting Wave Propagation Algorithm (DSWPA)**, where given the linear system of two-dimensional constant-coefficients

$$\mathbf{q}_t + \mathbf{A}\mathbf{q}_x + \mathbf{B}\mathbf{q}_y = 0, \tag{21}$$

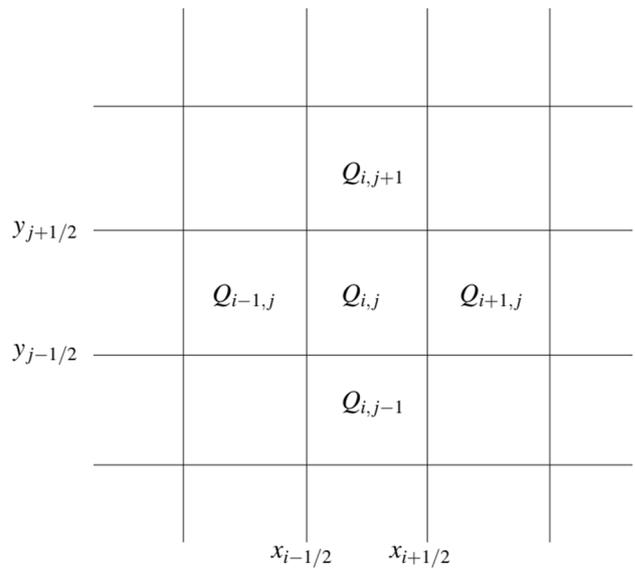

**Fig. 2** Two-dimensional finite volume grid, where $Q_{i,j}$ represents the average of cells

can be divided into two steps

$$\begin{aligned}\text{Step 1 } (x-\text{sweeps}) &: \mathbf{q}_t + \mathbf{A}\mathbf{q}_x = 0,\\ \text{Step 2 } (y-\text{sweeps}) &: \mathbf{q}_t + \mathbf{B}\mathbf{q}_y = 0.\end{aligned} \quad (22)$$

In the $x-$sweeps we have

$$\mathbf{Q}^*_{i,j} = \mathbf{Q}_{i,j} - \frac{\Delta t}{\Delta x}(\mathbf{F}^n_{i+1/2,j} - \mathbf{F}^n_{i-2/2,j}), \quad (23)$$

where $\mathbf{F}^n_{i-1/2,j}$ is the WPA numeric flux for the one-dimensional problem between the cells $\mathscr{C}_{i-1,j}$ and $\mathscr{C}_{ij}$ and in the $y-$sweeps we have

$$\mathbf{Q}^{n+1}_{i,j} = \mathbf{Q}^*_{i,j} - \frac{\Delta t}{\Delta y}(\mathbf{G}^*_{i,j+1/2} - \mathbf{G}^*_{i,j-1/2}), \quad (24)$$

where $\mathbf{G}^n_{i,j-1/2}$ is the WPA numeric flux for the one-dimensional problem between the cells $\mathscr{C}_{i,j-1}$ and $\mathscr{C}_{ij}$.

## 3.2 Semi-discrete Central-Upwind Scheme with Kurganov-Lin flux (CUp)

In this section, we describe the semi-discrete Central-Upwind scheme with Kurganov-Lin flux, developed by [12], which is based on the REA algorithm.

For simplicity, a uniform grid is considered, $x_i := i\Delta x$. It is also assumed that at the current time level, $t = t_n$, the average solution cells, $\mathbf{Q}^n_i$, are available. Then, the evolution of the solution computed for the next time level, $t = t_{n+1}$, can be presented by the one-dimensional semi-discrete scheme **CUp** that has the form of flux difference

$$\frac{d}{dt}\mathbf{Q}_i(t) = -\frac{\mathbf{H}_{i+1/2}(t) - \mathbf{H}_{i-1/2}(t)}{\Delta x}, \quad (25)$$

choosing the numerical flux $\mathbf{H}_{i-1/2}$, as

$$\mathbf{H}_{i-1/2}(t) := \frac{a^+_{i-1/2}\mathbf{f}(\mathbf{q}^-_{i-1/2}) - a^-_{i-1/2}\mathbf{f}(\mathbf{q}^+_{i-1/2})}{a^+_{i-1/2} - a^-_{i-1/2}} + a^+_{i-1/2}a^-_{i-1/2}\left[\frac{\mathbf{q}^+_{i-1/2} - \mathbf{q}^-_{i-1/2}}{a^+_{i-1/2} - a^-_{i-1/2}} - \mathbf{q}^{int}_{i-1/2}\right], \quad (26)$$

we obtain the Central-Upwind Kurganov-Lin scheme. Where the one-sided local velocities, $a^{\pm}_{i-1/2}$, are given by

$$a^+_{i-1/2} := \max_{w \in C(\mathbf{q}^-_{i-1/2}, \mathbf{q}^+_{i-1/2})} \left\{ \lambda_m \left( \frac{\partial \mathbf{f}}{\partial \mathbf{q}}(w) \right), 0 \right\} \geq 0, \tag{27}$$

$$a^-_{i-1/2} := \min_{w \in C(\mathbf{q}^-_{i-1/2}, \mathbf{q}^+_{i-1/2})} \left\{ \lambda_1 \left( \frac{\partial \mathbf{f}}{\partial \mathbf{q}}(w) \right), 0 \right\} \leq 0, \tag{28}$$

where $\lambda_1 < \cdots \lambda_m$ are the eigenvalues of the Jacobian $\frac{\partial \mathbf{f}}{\partial \mathbf{q}}$, and $C(\mathbf{q}^-_{i-1/2}, \mathbf{q}^+_{i-1/2})$ is a curve in phase space connecting the corresponding values on the left and right sides of the linear reconstruction at $x = x_{i-1/2}$,

$$\mathbf{q}^-_{i-1/2} := \mathbf{Q}_i - \frac{\Delta x}{2}(\mathbf{q}_x)^n_i, \tag{29}$$

$$\mathbf{q}^+_{i-1/2} := \mathbf{Q}_{i-1} + \frac{\Delta x}{2}(\mathbf{q}_x)^n_{i-1}, \tag{30}$$

respectively, where the numerical derivatives

$$\partial \mathbf{q}^n_i = \left. \frac{\partial \mathbf{q}}{\partial x} \right|_{x=x_i, t=t_n} + \mathcal{O}(\Delta x), \tag{31}$$

are calculated using flux limiters.

The flux limiter used in this work was the SuperBee (SB) defined by

$$\partial \mathbf{q}^n_i = maxmod \left( minmod \left( 2\frac{\mathbf{Q}^n_i - \mathbf{Q}^n_{i-1}}{\Delta x}, \frac{\mathbf{Q}^n_{i+1} - \mathbf{Q}^n_i}{\Delta x} \right), \right.$$
$$\left. minmod \left( \frac{\mathbf{Q}^n_i - \mathbf{Q}^n_{i-1}}{\Delta x}, 2\frac{\mathbf{Q}^n_{i+1} - \mathbf{Q}^n_i}{\Delta x} \right) \right), \tag{32}$$

where

$$maxmod(x_1, x_2, \cdots) = \begin{cases} \max_i\{x_i\}, & \text{if } x_i > 0 \; \forall i, \\ \min_i\{x_i\}, & \text{if } x_i < 0 \; \forall i, \\ 0, & \text{other case} \end{cases} \tag{33}$$

Finally, the correction term in (26) (which is an embedded anti-diffusion term that corresponds to the reduced value in numerical dissipation compared to the original semi-discrete CUp scheme of [24]) is

$$\mathbf{q}^{int}_{i-1/2} = minmod \left( \frac{\mathbf{q}^+_{i-1/2} - \mathbf{w}^{int}_{i-1/2}}{a^+_{i-1/2} - a^-_{i-1/2}}, \frac{\mathbf{w}^{int}_{i-1/2} - \mathbf{q}^-_{i-1/2}}{a^+_{i-1/2} - a^-_{i-1/2}} \right), \tag{34}$$

with

$$\mathbf{w}^{int}_{i-1/2} = \frac{a^+_{i-1/2}\mathbf{q}^+_{i-1/2} - a^-_{i-1/2}\mathbf{q}^-_{i-1/2} - \left\{ \mathbf{f}(\mathbf{q}^+_{i-1/2}) - \mathbf{f}(\mathbf{q}^-_{i-1/2}) \right\}}{a^+_{i-1/2} - a^-_{i-1/2}}. \tag{35}$$

For the two-dimensional Central-Upwind Semi-Discrete Scheme with Kurganov-Lin flux, similar to case 1-D, we consider a uniform grid $x_i = i\Delta x$, $y_j = j\Delta y$, and $\Delta t := t_{n+1} - t_n$. The resulting semi-discrete two-dimensional CUp scheme will then be obtained in the following flow difference form:

$$\frac{d}{dt}\mathbf{Q}_{i,j}(t) = -\frac{\mathbf{H}^x_{i+1/2,j}(t) - \mathbf{H}^x_{i-1/2,j}(t)}{\Delta x} - \frac{\mathbf{H}^y_{i,j+1/2}(t) - \mathbf{H}^y_{i,j-1/2}(t)}{\Delta y}, \tag{36}$$

choosing second-order numeric flux as being

$$\mathbf{H}^x_{i-1/2,j}(t) := \frac{a^+_{i-1/2,j}\mathbf{f}(\mathbf{q}^E_{i-1,j}) - a^-_{i-1/2,j}\mathbf{f}(\mathbf{q}^W_{i,j})}{a^+_{i-1/2,j} - a^-_{i-1/2,j}} \\ + a^+_{i-1/2,j}a^-_{i-1/2,j}\left[\frac{\mathbf{q}^W_{i,j} - \mathbf{q}^E_{i-1,j}}{a^+_{i-1/2,j} - a^-_{i-1/2,j}} - \mathbf{q}^{x-int}_{i-1/2,j}\right], \quad (37)$$

$$\mathbf{H}^y_{i,j-1/2}(t) := \frac{b^+_{i,j-1/2}\mathbf{f}(\mathbf{q}^N_{i,j-1}) - b^-_{i,j-1/2}\mathbf{f}(\mathbf{q}^S_{i,j})}{b^+_{j-1/2} - b^-_{i,j-1/2}} \\ + b^+_{i,j-1/2}b^-_{i,j-1/2}\left[\frac{\mathbf{q}^S_{i,j} - \mathbf{q}^N_{i,j-1}}{b^+_{i,j-1/2} - b^-_{i,j-1/2}} - \mathbf{q}^{y-int}_{i,j-1/2}\right], \quad (38)$$

we get the CUp scheme with Kurganov-Lin flux, unilateral local propagation speeds can be estimated, for example, by

$$\begin{aligned}
a^+_{i-1/2,j} &:= \max\left\{\lambda_1\left(\mathbf{A}(\mathbf{q}^W_{i,j})\right), \lambda_m\left(\mathbf{A}(\mathbf{q}^E_{i-1,j})\right), 0\right\}, \\
b^+_{i,j-1/2} &:= \max\left\{\lambda_1\left(\mathbf{B}(\mathbf{q}^S_{i,j})\right), \lambda_m\left(\mathbf{B}(\mathbf{q}^N_{i-1,j})\right), 0\right\}, \\
a^-_{i-1/2,j} &:= \min\left\{\lambda_1\left(\mathbf{A}(\mathbf{q}^W_{i,j})\right), \lambda_m\left(\mathbf{A}(\mathbf{q}^E_{i-1,j})\right), 0\right\}, \\
b^-_{i-1/2,j} &:= \min\left\{\lambda_1\left(\mathbf{B}(\mathbf{q}^S_{i,j})\right), \lambda_m\left(\mathbf{B}(\mathbf{q}^N_{i-1,j})\right), 0\right\},
\end{aligned} \quad (39)$$

where, $\lambda_1 < \lambda_2 < \cdots < \lambda_m$ are the $m$ eigenvalues of the corresponding Jacobians, $\mathbf{A} := \frac{\partial \mathbf{f}}{\partial \mathbf{q}}$ and $\mathbf{B} := \frac{\partial \mathbf{g}}{\partial \mathbf{q}}$, and the point values of linear piecewise reconstruction are given by

$$\begin{aligned}
\mathbf{q}^E_{i,j} &:= \mathbf{Q}^n_{i,j} + \tfrac{\Delta x}{2}(\mathbf{q}_x)^n_{i,j} & \mathbf{q}^N_{i,j} &:= \mathbf{Q}^n_{i,j} + \tfrac{\Delta y}{2}(\mathbf{q}_y)^n_{i,j}, \\
\mathbf{q}^W_{i,j} &:= \mathbf{Q}^n_{i,j} - \tfrac{\Delta x}{2}(\mathbf{q}_x)^n_{i,j} & \mathbf{q}^S_{i,j} &:= \mathbf{Q}^n_{i,j} - \tfrac{\Delta y}{2}(\mathbf{q}_y)^n_{i,j},
\end{aligned} \quad (40)$$

where again the numerical derivatives are calculated using flux limiters.

Similarly to the one-dimensional case, we call $\mathbf{q}^{x-int}_{i-1/2,j}$ and $\mathbf{q}^{y-int}_{i,j-1/2}$ the dissipation correction terms that are given by

$$\mathbf{q}^{x-int}_{i-1/2,j} = \text{minmod}\Bigg(\frac{\mathbf{q}^{NW}_{i,j} - \mathbf{w}^{int}_{i-1/2,j}}{a^+_{i-1/2,j} - a^-_{i-1/2,j}}, \frac{\mathbf{w}^{int}_{i-1/2,j} - \mathbf{q}^{NE}_{i-1,j}}{a^+_{i-1/2,j} - a^-_{i-1/2,j}}, \\ \frac{\mathbf{q}^{SW}_{i,j} - \mathbf{w}^{int}_{i-1/2,j}}{a^+_{i-1/2,j} - a^-_{i-1/2,j}}, \frac{\mathbf{w}^{int}_{i-1/2,j} - \mathbf{q}^{SE}_{i-1,j}}{a^+_{i-1/2,j} - a^-_{i-1/2,j}}\Bigg), \quad (41)$$

$$\mathbf{q}^{y-int}_{i,j-1/2} = \text{minmod}\Bigg(\frac{\mathbf{q}^{SW}_{i,j} - \mathbf{w}^{int}_{i,j-1/2}}{a^+_{i,j-1/2} - a^-_{i,j-1/2}}, \frac{\mathbf{w}^{int}_{i,j-1/2} - \mathbf{q}^{NW}_{i,j-1}}{a^+_{i-1/2,j} - a^-_{i,j-1/2}}, \\ \frac{\mathbf{q}^{SE}_{i,j} - \mathbf{w}^{int}_{i,j-1/2}}{a^+_{i,j-1/2} - a^-_{i,j-1/2}}, \frac{\mathbf{w}^{int}_{i,j-1/2} - \mathbf{q}^{NE}_{i,j-1}}{a^+_{i,j-1/2} - a^-_{i,j-1/2}}\Bigg), \quad (42)$$

where

$$\begin{aligned}
\mathbf{w}^{int}_{i-1/2,j} &= \frac{a^+_{i-1/2,j}\mathbf{q}^W_{i,j} - a^-_{i-1/2,j}\mathbf{q}^E_{i-1,j} - [\mathbf{f}(\mathbf{q}^W_{i,j}) - \mathbf{f}(\mathbf{q}^E_{i-1,j})]}{a^+_{i-1/2,j} - a^-_{i-1/2,j}}, \\
\mathbf{w}^{int}_{i,j-1/2} &= \frac{b^+_{i,j-1/2}\mathbf{q}^S_{i,j} - b^-_{i,j-1/2}\mathbf{q}^N_{i,j-1} - [\mathbf{f}(\mathbf{q}^S_{i,j}) - \mathbf{f}(\mathbf{q}^N_{i,j-1})]}{b^+_{i,j-1/2} - b^-_{i,j-1/2}},
\end{aligned} \quad (43)$$

and $\mathbf{q}_{i,j}^{NE}, \mathbf{q}_{i,j}^{NW}, \mathbf{q}_{i,j}^{SE}, \mathbf{q}_{i,j}^{SW}$ are the values of the corresponding corner points of the linear reconstruction by parts in cell $(i, j)$ given by

$$\begin{aligned}
\mathbf{q}_{i,j}^{NE} &:= \mathbf{Q}_{i,j}^n + \frac{\Delta x}{2}(\mathbf{q}_x)_{i,j}^n + \frac{\Delta y}{2}(\mathbf{q}_y)_{i,j}^n, & \mathbf{q}_{i,j}^{SE} &:= \mathbf{Q}_{i,j}^n + \frac{\Delta x}{2}(\mathbf{q}_x)_{i,j}^n - \frac{\Delta y}{2}(\mathbf{q}_y)_{i,j}^n, \\
\mathbf{q}_{i,j}^{NW} &:= \mathbf{Q}_{i,j}^n - \frac{\Delta x}{2}(\mathbf{q}_x)_{i,j}^n + \frac{\Delta y}{2}(\mathbf{q}_x)_{i,j}^n, & \mathbf{q}_{i,j}^{SW} &:= \mathbf{Q}_{i,j}^n - \frac{\Delta x}{2}(\mathbf{q}_x)_{i,j}^n - \frac{\Delta y}{2}(\mathbf{q}_y)_{i,j}^n.
\end{aligned} \quad (44)$$

## 4 Numerical results

In this section, we introduce the numerical seismic test cases utilized to assess the performance and competitiveness of the Godunov-type finite volume methods presented earlier, particularly in comparison with the staggered finite difference method. Throughout these examples, we employ reflective or solid wall boundary conditions at the top and bottom and extend the domain on both sides to prevent reflections at these boundaries from interfering with the receivers. This approach allows us to avoid the use of absorbing boundary conditions, which are not the primary focus of our study. All simulations were performed using a numerical server with 2x Intel Xeon X5690 (6c/12t) 3.47 GHz processors and 64 GB RAM.

### 4.1 Test case 1—Heterogeneous parallel velocity model (Synthetic test case)

In the first test case, our domain spans from $x_0 = 0$ to $x_1 = 1000$ m and $y_0 = 0$ to $y_1 = 1000$ m, featuring a velocity profile as depicted in Fig. 3. This profile ranges from velocities of 1.5 km/s to 5.5 km/s. We utilize a Ricker wavelet given by:

$$g(t) = (1 - 2\pi^2 f_M^2 (t - t_0)^2) \exp(\pi^2 f_M^2 (t - t_0)^2),$$

where $f_M = 0.015$ kHz represents the peak frequency and $t_0 = 1/f_M$. The seismic source is positioned at the domain's center, at a depth of 20 m, while receivers are placed along the x-axis at a depth of 20 m.

In Fig. 3, we display the velocity profile, highlighting the positions of the seismic source (red marker) and the receivers (green markers), which are spaced at intervals of 50 m. Figure 4 illustrates the reference solution for the field $\sigma(x, y)$ and the associated reference seismogram. The subfigures show a reference snapshot of the wavefield at $t = 0.3$ s (panel A), and a 1D vertical slice of the wavefield at $x = 330$ m (panel B). The reference seismogram at $t = 1$ s is shown in panel (C) alongside with the seismogram data recorded at a receiver positioned at $x = 330$ m and $y = 20$ m (panel D). This solution was generated using a 20th-order finite difference method implemented with the Devito package [25, 26]. The spatial

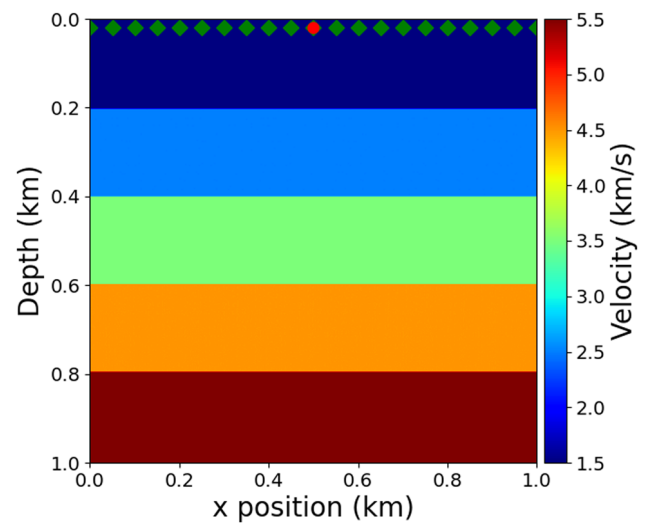

**Fig. 3** **Test case 1.** Heterogeneous parallel velocity model. The domain spans from $x_0 = 0$ to $x_1 = 1000$ m and $y_0 = 0$ to $y_1 = 3000$ m

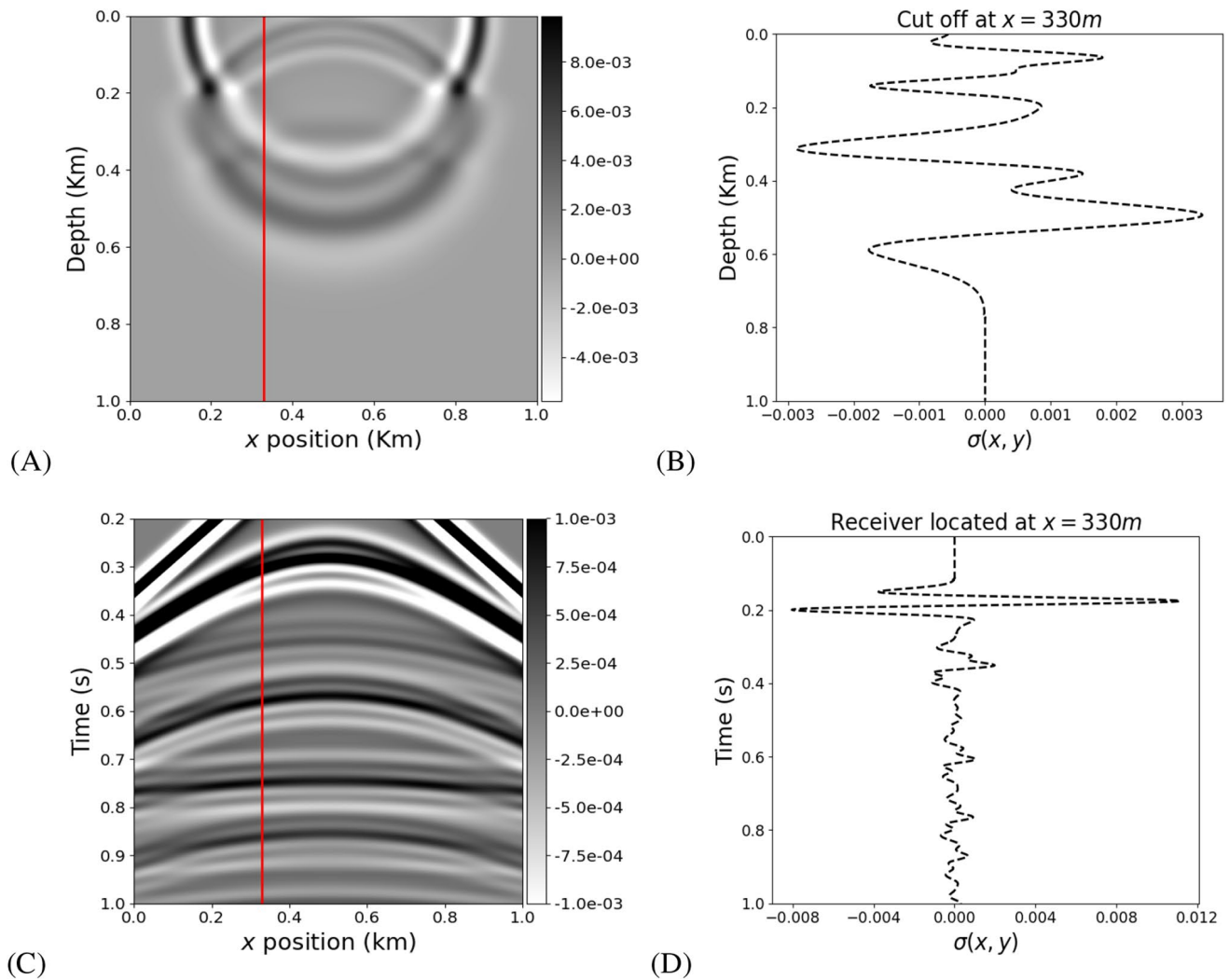

**Fig. 4** **Test case 1: Reference.** Top row: **A** Reference snapshot of the wavefield at time $t = 0.3$ s. **B** 1D vertical snapshot (cut) of the wavefield at $x = 330$ m from **A**, corresponding to the red line indicated in the snapshot. Bottom row: **C** Reference seismogram up until time $t = 1$ s. **D** The Seismogram data recorded at the receiver located at $x = 330$ m and $y = 20$ m. This solution was computed using a 20th-order finite difference method implemented with the Devito package, with a spatial resolution $\Delta x = \Delta y = 0.15625$ m. The Courant-Friedrichs-Lewy (CFL) number was set to 0.25

resolution was set to $\Delta x = \Delta y = 0.15625$ m, and the Courant-Friedrichs-Lewy (CFL) number was 0.25. This time and spatial resolution, along with the high-order scheme, generates a sufficiently accurate solution to use as a reference.

Figure 5 shows the absolute differences between the reference wavefield snapshot (shown in Fig. 4A) and the snapshots obtained by the DF2, DF8, DSWPA, FWPA and CUP methods with a grid spacing of $\Delta x = \Delta y = 5$ m. In the top row of Fig. 6, we compare these methods to the reference solution along the vertical cut at $x = 330$ m (shown in Fig. 4A). According to these figures, we can see that the finite difference methods DF2 and DF8 exhibit slight delays relative to the reference solution and the finite volume methods, suggesting a slightly higher dispersion. However, we can notice that the finite volume method shows relatively early arrivals with opposite first motions, also indicative of dispersion although smaller than the finite difference methods.

On the other hand, the finite volume methods show more diffusion. This is further illustrated in the bottom row of Fig. 6, where the max-norm errors are consistently higher for the finite volume methods, while the 1-norm errors are lower. The plot also includes slope 1 and slope 2 lines, representing the 1st and 2nd order of convergence, respectively, providing a reference for evaluating the error decay rates of the different methods and we can see that all the methods show a first order convergence although the theoretical order of each method is higher, this is due to the heterogeneity of the velocity profile and the use of the Total Variation Diminishing (TVD) limiters in the finite volume methods.

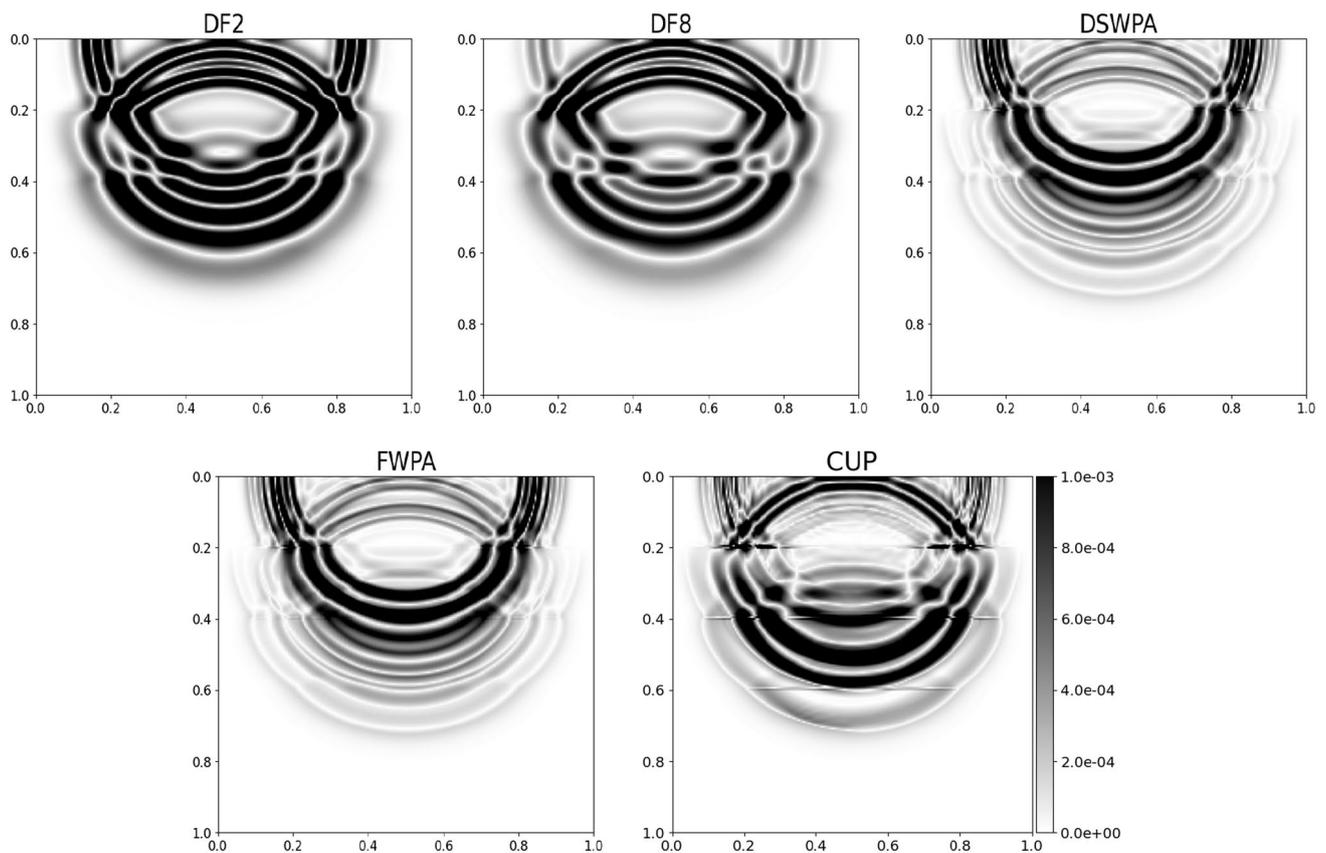

**Fig. 5** **Test case 1: Error snapshots.** Absolute value of the difference between the reference solution and the solution obtained by implementing the DF2, DF8, DSWPA, FWPA, and CUP methods with $\Delta x = \Delta y = 5$ m

Similarly, Fig. 7 depicts the absolute differences between the reference seismogram and those obtained using the DF2, DF8, DSWPA, FWPA, and CUP methods for the same grid spacing at $t = 1$ s. Figure 8 shows the comparison with the reference seismogram recorded at the receiver located at $x = 330$ m (shown on the bottom right panel of Fig. 4). Here, the dispersion introduced by DF2 and DF8 is more apparent, and the higher diffusivity of finite volume methods is highlighted compared to finite difference methods.

Furthermore, Fig. 9 (left) illustrates the computational time required by each method to obtain the seismogram for test case 1. As expected, finite volume methods are more computationally expensive than finite difference methods. Specifically, the computational time required by finite volume methods DFWPA and FWPA on a grid with $\Delta x = \Delta y$ is approximately equivalent to that required by the finite difference method of order 8 on a grid with $\Delta x/2 = \Delta y/2$. Figure 9 (A and B) shows the error obtained by each method in 1-norm and max-norm versus the computational time required. It is evident that in 1-norm, the time required by the DF8 method to achieve a certain error is comparable to that required by DSWPA and FWPA for test case 1. However, as expected, in max-norm, the finite volume methods show inferior performance.

In Figs. 10 and 11 we introduced the parameter $\lambda$, representing the percentage decrease in discontinuity in velocity at each layer of the velocity profile presented in this test case, in these figures we show the relationship between $\lambda$ (where $\lambda = 0\%$ corresponds to the original velocity profile with increments of 1 km/s at each interface) and the error of each method in both the 1-norm and the max-norm for fixed $\Delta x$ and $\Delta y$. We observe that as $\lambda$ increases, indicating a reduction in the discontinuity at each interface, the error obtained with the finite difference method is significantly diminished, this shows that the greater the jump in the discontinuity in the velocity profile, the greater the error introduced by finite difference method.

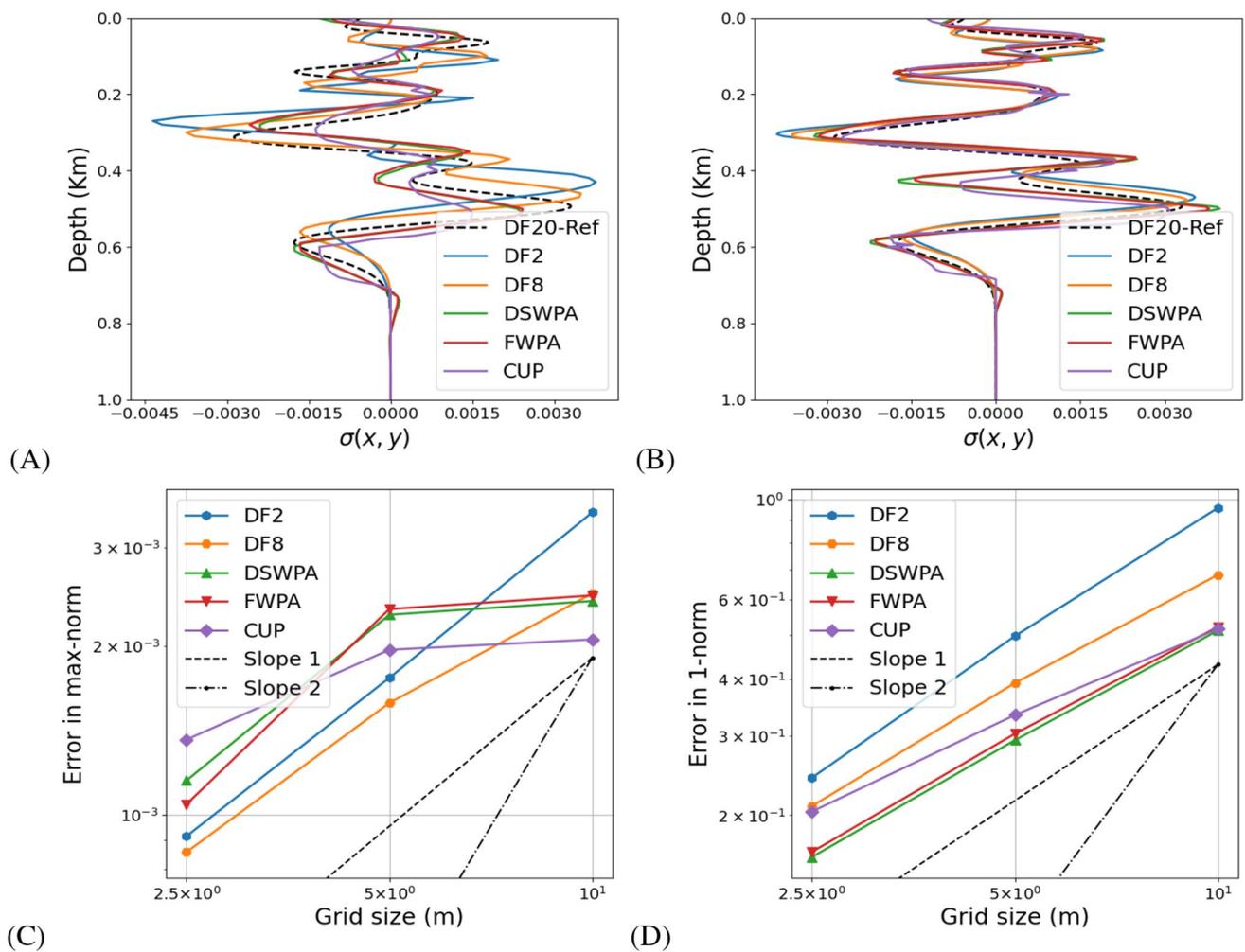

**Fig. 6 Test case 1.** Top row: a comparison with the reference solution (snapshot of the wavefield at $t = 0.3$ s along the vertical cut at $x = 330$ m) is presented for the DF2, DF8, DSWPA, FWPA, and CUP methods using $\Delta x = \Delta y = 5$ m (**A**) and $\Delta x = \Delta y = 2.5$ m (**B**). Bottom row: log-log plots depict the error versus the grid size for the DF2, DF8, DSWPA, FWPA, and CUP methods using the max-norm (**C**) and 1-norm (**D**). Here slope 1 and slope 2 lines represent references for 1st and 2nd order of convergence, respectively

### 4.2 Test case 2—SEG/EAGE salt body velocity profile.

In this test case, the domain spans from $x_0 = 2000$ m to $x_1 = 11,000$ m and $y_0 = 0$ to $y_1 = 3000$ m, with a velocity profile depicted in Fig. 12, with velocities ranging from 1.5 km/s to 4.5 km/s.

The source was positioned at the middle of the domain at a depth of 40 m, while receivers were situated along the $x$ axis at a depth of 40 m and a frequency $f_M = 0.015$ KHz.

Figure 13 displays the absolute difference between the reference seismogram for test case 2 and the seismogram obtained by DF2, DF8, DSWPA, FWPA, and CUP methods with $\Delta x = \Delta y = 20$ m at $t = 2$ s. Figure 14 shows the comparison with the reference seismogram recorded at the receiver located at $x = 7000$ m and $t = 2$ s in the top row. In the bottom row, on the left and center, we show the log-log plots of the error versus grid size, and on the right side of the figure, we display a plot of computational time versus grid size. The reference solution was computed using a 20th-order finite difference method with a spatial resolution of $\Delta x = \Delta y = 1.25$ m.

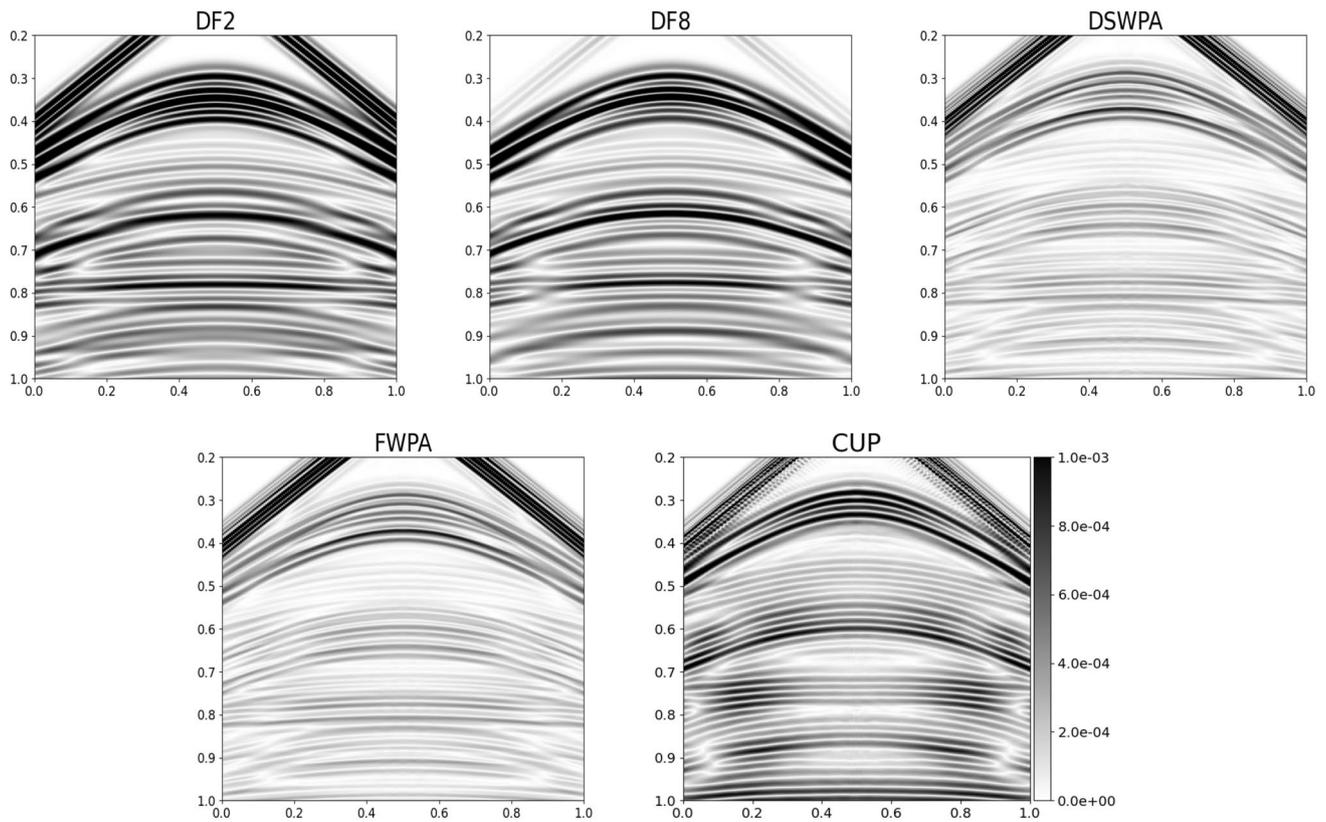

**Fig. 7** **Test case 1.** Absolute difference between the reference seismogram and the seismograms obtained using the DF2, DF8, DSWPA, FWPA, and CUP methods, using a spatial resolution of $\Delta x = \Delta y = 5$ m

### 4.3 Test case 3—Marmousi velocity profile

In this test case, we utilized the Marmousi velocity profile, depicted in Fig. 15, originally spanning a horizontal extension of 17 km and a depth of 3.5 km. Simulations were conducted within a subdomain ranging from $x_0 = 4000$ m to $x_1 = 14{,}000$ m, maintaining the same depth as the original profile ($y_0 = 0$ to $y_1 = 3500$ m). The source was positioned at the middle of the domain at a depth of 40 m, while receivers were distributed along the $x$ axis at the same depth. A frequency of $f_M = 0.005$ KHz was employed.

Figure 16 exhibits the absolute difference between the reference seismogram for test 3 and the seismogram obtained using the DF2, DF8, DSWPA, FWPA, and CUP methods with $\Delta x = \Delta y = 20$ m at $t = 2.5$ s. Figure 14 shows the comparison with the reference solution calculated with a finite difference method of order 20 using the receiver located at $x = 7000$ m at the top of the figure, at left and center we show the log-log plots of the error versus the grid size and in the right of the figure, we show a plot of the computational time versus grid size. The reference solution was computed with a 20th-order finite difference method with a spatial resolution of $\Delta x = \Delta y = 1.25$ m.

### 4.4 Test case 4—a typical velocity field of Santos Basin

For the fourth test, we utilized a typical velocity profile from the Santos Basin (as depicted in Fig. 18), an expansive oil field covering 352,260 square kilometers situated in the Atlantic Ocean, approximately 300 km southeast of Santos, Brazil. Renowned as one of the largest ocean basins in the country, the Santos Basin harbors numerous oil reserve exploration and exploitation sites.

The velocity profile initially spanned a horizontal distance of 70.3 km with a depth of 9.92 km. However, our simulations focused on a subdomain within $x = [20{,}000$ m$, 50{,}000$ m$]$, covering a horizontal extent of 30 km, and $y = [0, 8000$ m$]$, with a depth of 8 km. The seismic source was positioned at the midpoint of the subdomain, specifically at $x = 35{,}000$ m and a depth of 32 m. Receivers were aligned along the $x$ axis at the same depth. We employed a frequency of $f_M = 0.005$ KHz.

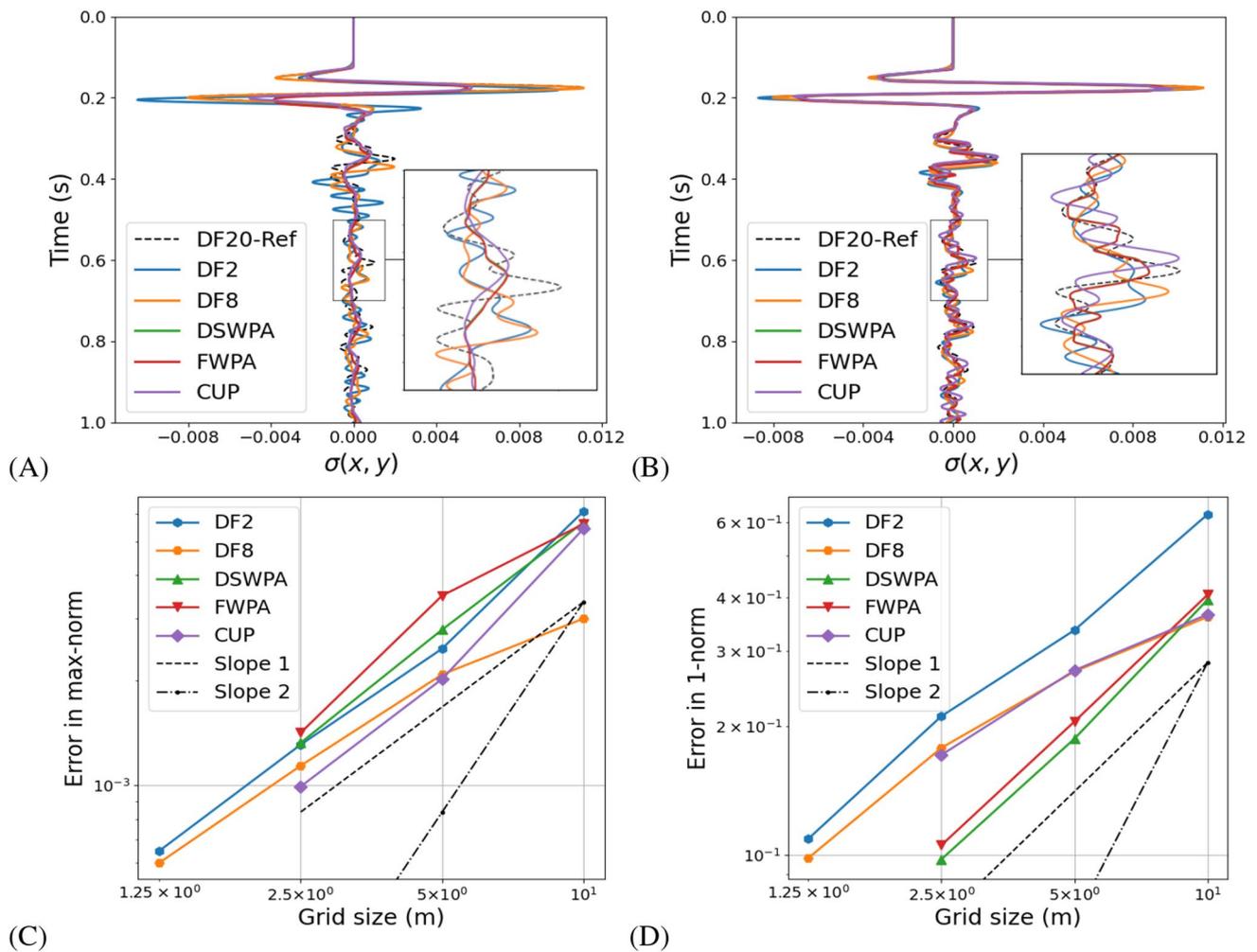

**Fig. 8 Test case 1.** Top row: a comparison with the reference seismogram recorded at the receiver located at $x = 330$ m and $t = 1$ s is presented for the DF2, DF8, DSWPA, FWPA, and CUP methods with $\Delta x = \Delta y = 5$ m (**A**) and $\Delta x = \Delta y = 2.5$ m (**B**). Bottom row: log-log plots depict the error versus the grid size for the same methods using the max-norm (**C**) and 1-norm (**D**). Here slope 1 and slope 2 lines represent references for 1st and 2nd order of convergence, respectively

Figure 19 showcases a reference seismogram at $t = 5$ s, calculated with a 20th-order finite difference method in a grid with 9601 points in $x$ and 4001 points in $y$. Additionally, Fig. 19 presents the reference seismogram recorded for the receiver located at $x = 33{,}350$ m at $t = 5$ s. Figure 20 illustrates the absolute difference between the reference seismogram for test case 4 and the seismogram obtained using the DF2, DF8, DSWPA, FWPA, and CUP methods with $\Delta x = 25$ m and $\Delta y = 8$ m at $t = 5$ s. Additionally, it presents a comparison with the reference solution using the receiver located at $x = 33{,}350$ m.

We can conclude from the Figs. 14, 17 and 20 that when implementing finite volume methods in more realistic scenarios (test cases 2, 3, and 4), the finite difference methods DF2 and DF8 outperform the finite volume methods. Moreover, the finite volume methods require more computational time. These figures highlight that the finite volume methods exhibit greater error compared to the finite difference methods, both in the max-norm and the 1-norm, and necessitate more computational time, aligning with our expectations.

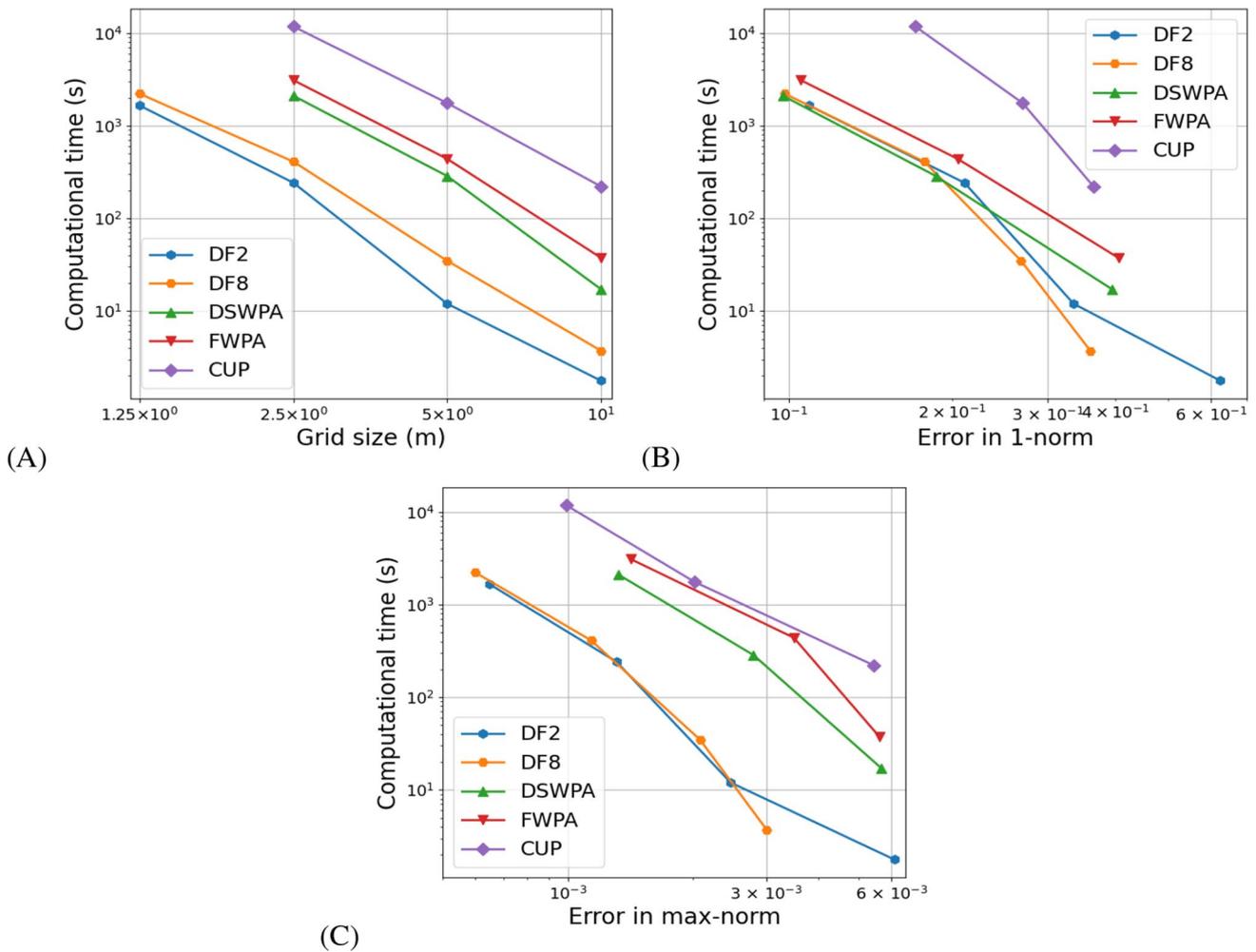

**Fig. 9** **Test case 1.** **A** plot of the computational time (in seconds) vs. grid size, **B** computational time vs. error in 1-norm for each grid size, **C** computational time vs. error in max-norm for each grid size, for the methods DF2, DF8, DSWPA, FWPA and CUP, final time $t = 1$ s

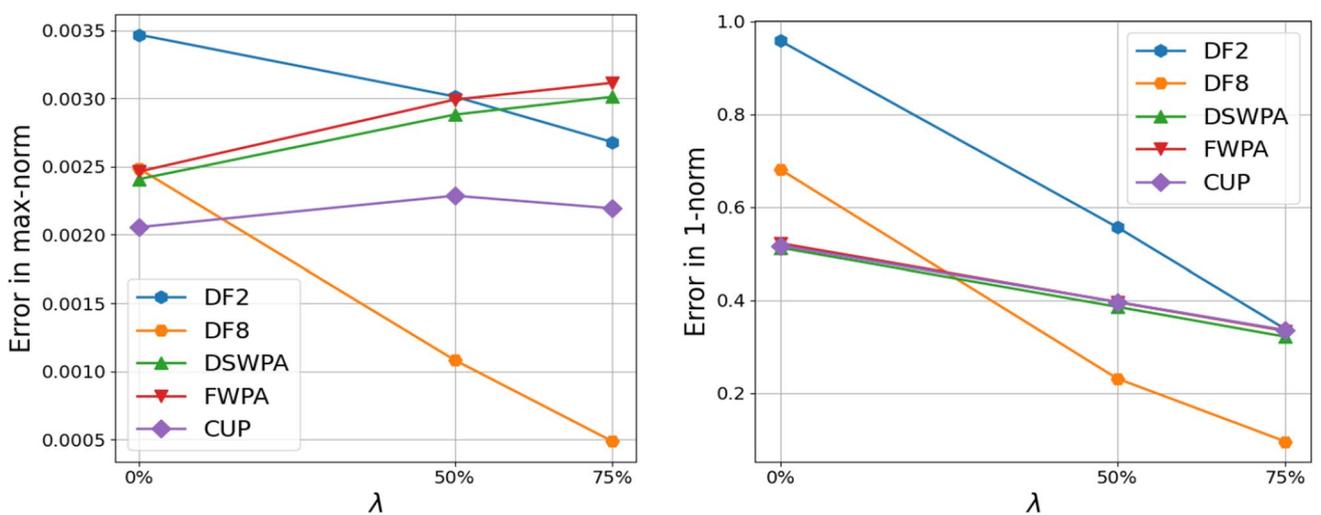

**Fig. 10** **Test case 1 (variation in discontinuity).** Plot of the error vs. the percentage of decrease in the discontinuity in velocity in each layer from the original velocity profile ($\lambda$) for the DF2, DF8, DSWPA, FWPA, and CUP methods. with a cutoff at $x = 330$ m at $t = 0.3$ s, the Grid spacing is $\Delta x = \Delta y = 10$ m, max-norm (left), 1-norm (right)

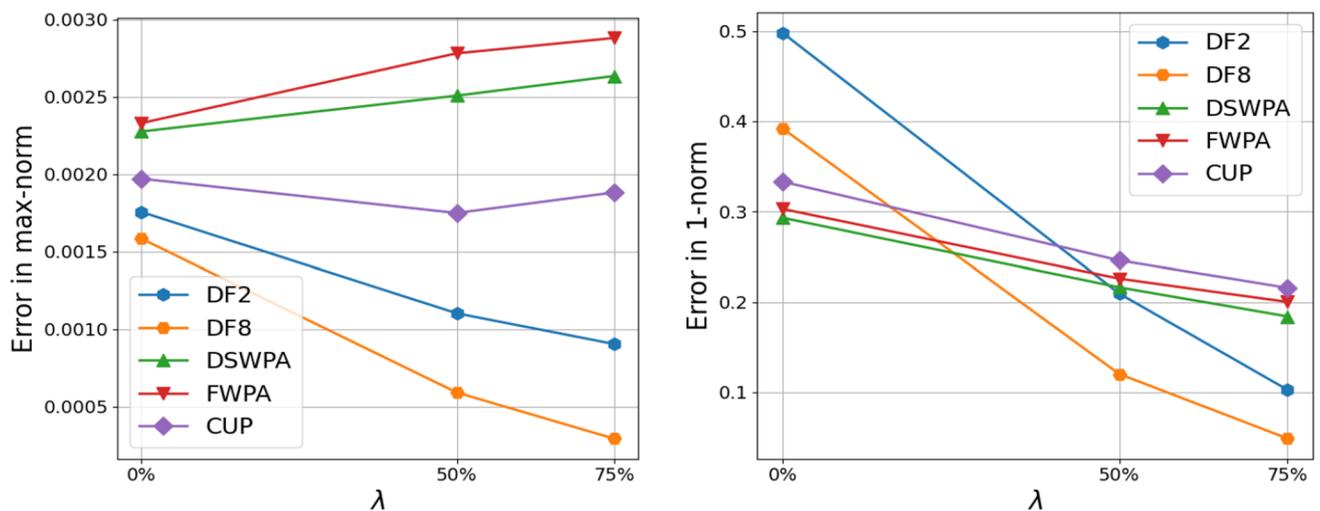

**Fig. 11** **Test case 1 (variation in discontinuity).** Plot of the error vs. the percentage of decrease in the discontinuity in velocity in each layer from the original velocity profile ($\lambda$) for the DF2, DF8, DSWPA, FWPA, and CUP methods. with a cutoff at $x = 330$ m at $t = 0.3$ s, the Grid spacing is $\Delta x = \Delta y = 5$ m, max-norm (left), 1-norm (right)

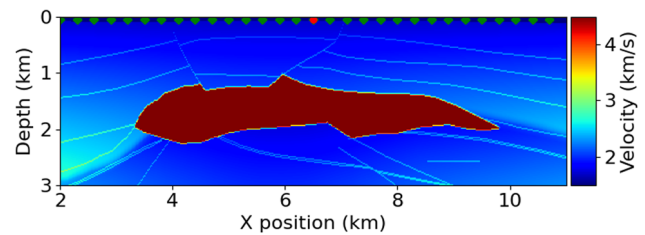

**Fig. 12** **Test case 2.** SEG/EAGE salt body velocity profile. The domain spans from $x_0 = 2000$ m to $x_1 = 11{,}000$ m and $y_0 = 0$ to $y_1 = 3000$ m

## 5 Conclusions

In our study, we conduct a comparative analysis of finite volume and finite difference methods applied to seismic wave propagation. This investigation addresses a notable gap in the literature, as finite volume methods, while extensively validated in fluid dynamics and hyperbolic conservation laws, remain underexplored in the context of seismic problems. By juxtaposing these methods with the widely-used finite difference techniques, we aim to highlight their respective strengths and limitations under varying seismic conditions.

Our findings indicate that the performance of these methods is highly dependent on the nature of the velocity profiles. In scenarios characterized by sharp velocity discontinuities at layer interfaces, finite volume methods demonstrate robustness and achieve lower numerical errors in 1-norm, despite their inherently higher computational overhead and increased numerical dissipation. Conversely, when velocity profile discontinuities are smoother or less pronounced, finite difference methods outperform finite volume approaches, offering superior computational efficiency and numerical accuracy, particularly in more realistic seismic simulations.

The finite volume methods analyzed in this study are high-resolution techniques, achieving up to second-order accuracy. This strikes a balance between computational cost and numerical precision, making them effective for complex heterogeneous media and intricate subsurface structures. Nevertheless, these methods can be extended to higher orders if necessary. For example, the Wave Propagation Algorithm (WPA) can achieve higher-order accuracy, as demonstrated in [27], and similar advancements have been made for the Central-Upwind method in [7]. These extensions enhance accuracy but come with a substantial increase in computational cost.

While finite volume methods perform well in capturing discontinuities and are well-suited for hyperbolic problems involving shock waves or rarefaction waves, their performance in linear seismic wave propagation scenarios is notably less efficient. This shortcoming stems from their higher computational cost, the challenges associated with their implementation, and the numerical dissipation introduced by TVD limiters. We note, however, that this comparison is inherently

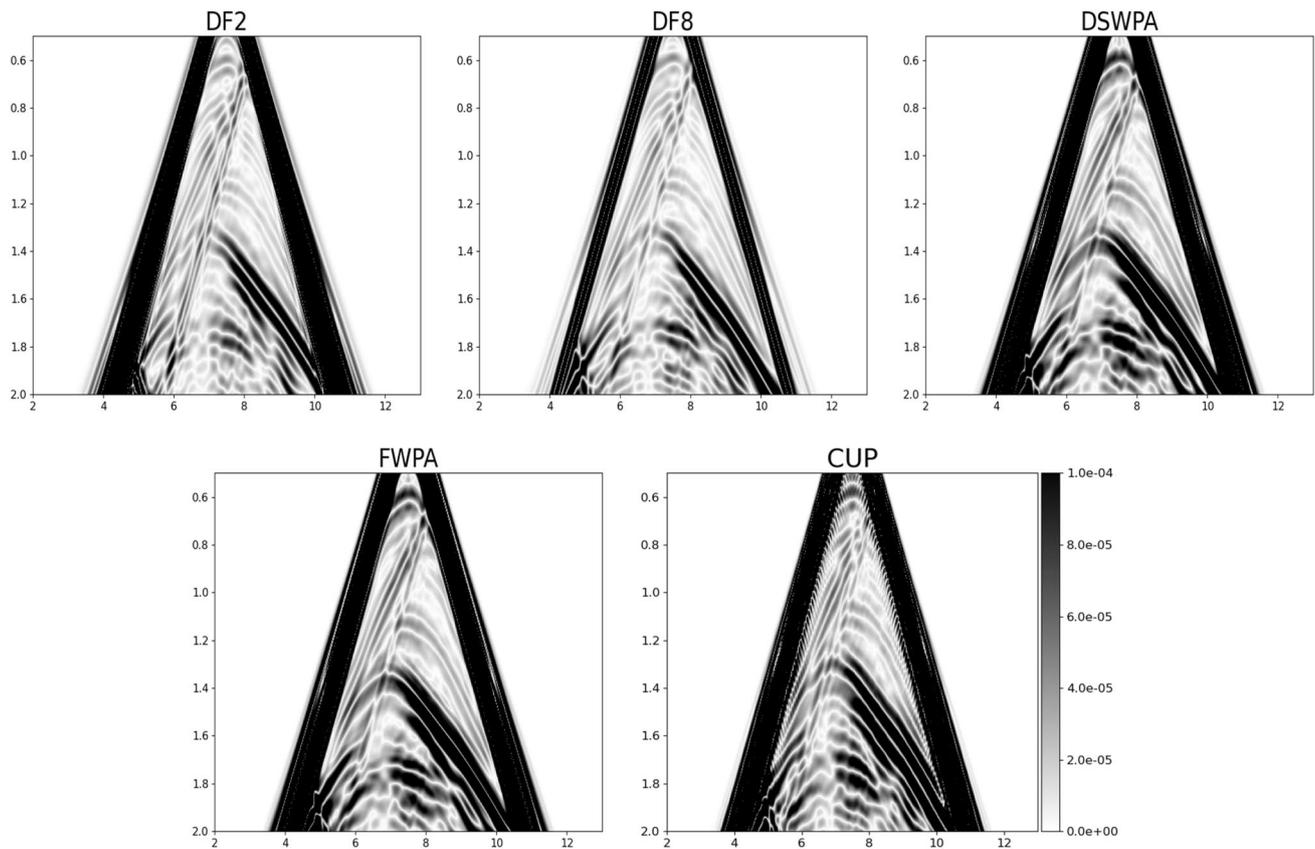

**Fig. 13 Test case 2 (SEG/EAGE).** Absolute difference between the reference seismogram and the seismograms obtained using the DF2, DF8, DSWPA, FWPA, and CUP methods. These computations were conducted with a spatial resolution of $\Delta x = \Delta y = 20$ m and a final time $t = 2$ s

tied to our specific implementations of each method, both of which were developed in Python using the NumPy library to maximize efficiency. Computational times were measured by averaging the runtimes of five simulations executed on the same machine, ensuring a consistent and fair basis for comparison.

In conclusion, finite difference methods remain the preferred choice for large-scale seismic wave propagation problems, offering an optimal balance between accuracy, computational efficiency, and ease of implementation. On the other hand, finite volume methods demonstrate clear advantages in scenarios with sharp discontinuities and highly heterogeneous media, yet their higher computational cost and dissipation issues limit their broader applicability in seismic modeling. Future research could focus on optimizing finite volume methods for linear wave propagation problems, exploring higher-order extensions to reduce numerical dissipation

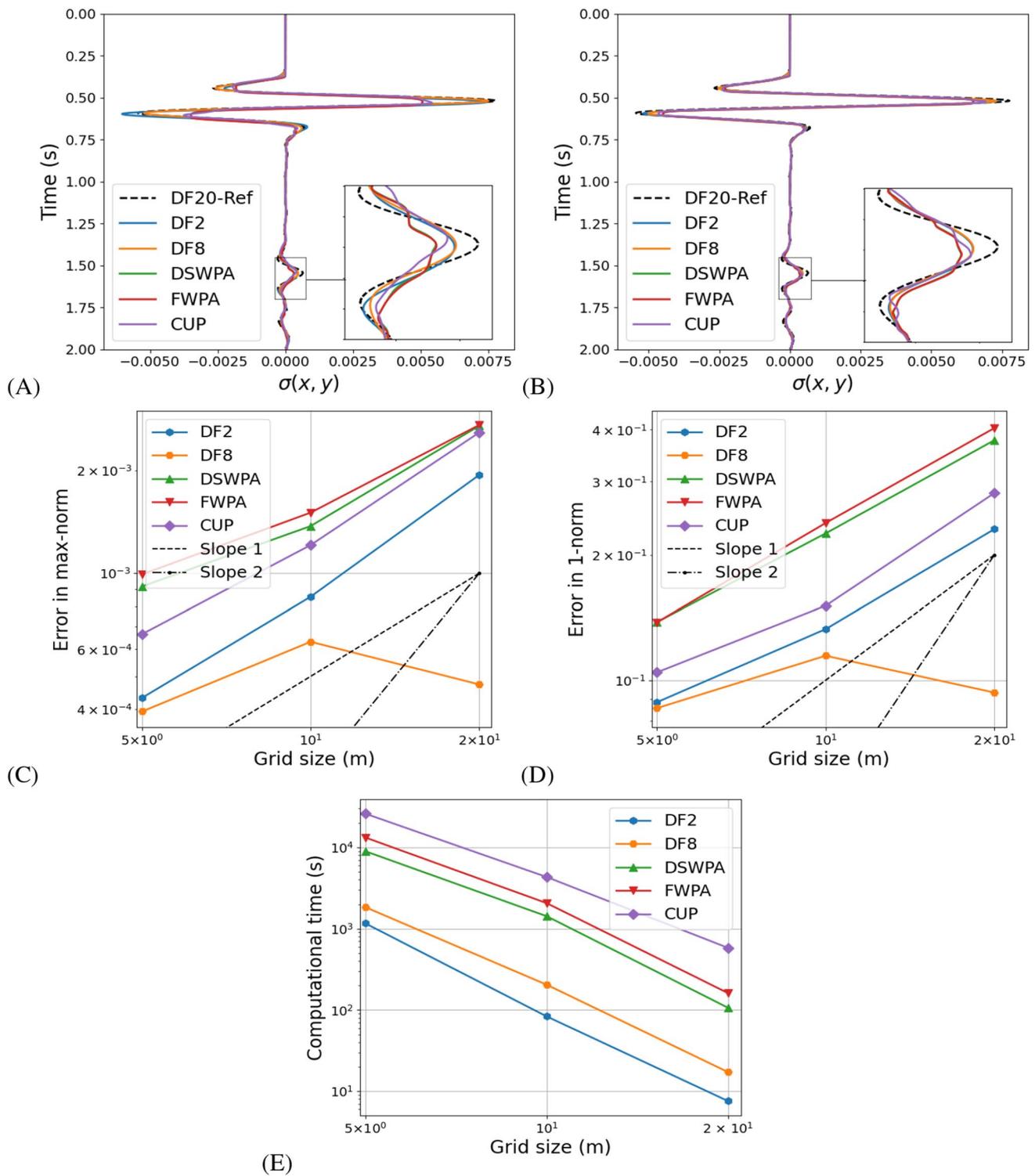

**Fig. 14 Test case 2 (SEG/EAGE).** Top row: we present a comparison with the reference seismogram recorded at the receiver located at $x = 7000$ m and $t = 2$ s for the DF2, DF8, DSWPA, FWPA, and CUP methods with $\Delta x = \Delta y = 20$ m (**A**) and $\Delta x = \Delta y = 10$ m (**B**). Bottom row: log-log plots show the error versus grid size for the same methods using the max-norm (**C**), 1-norm (**D**), and a plot of grid size versus computational time for this test case (**E**). The lines labeled slope 1 and slope 2 represent references for 1st and 2nd order of convergence, respectively

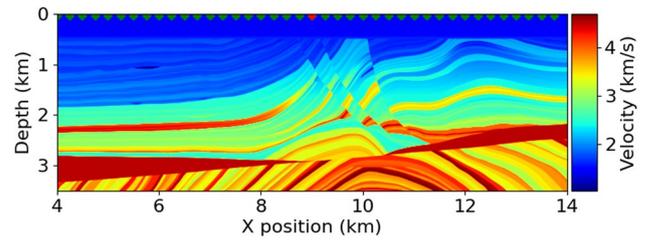

**Fig. 15** **Test case 3.** Marmousi velocity profile. The domain spans from $x_0 = 4000$ m to $x_1 = 14{,}000$ m and $y_0 = 0$ to $y_1 = 3500$ m

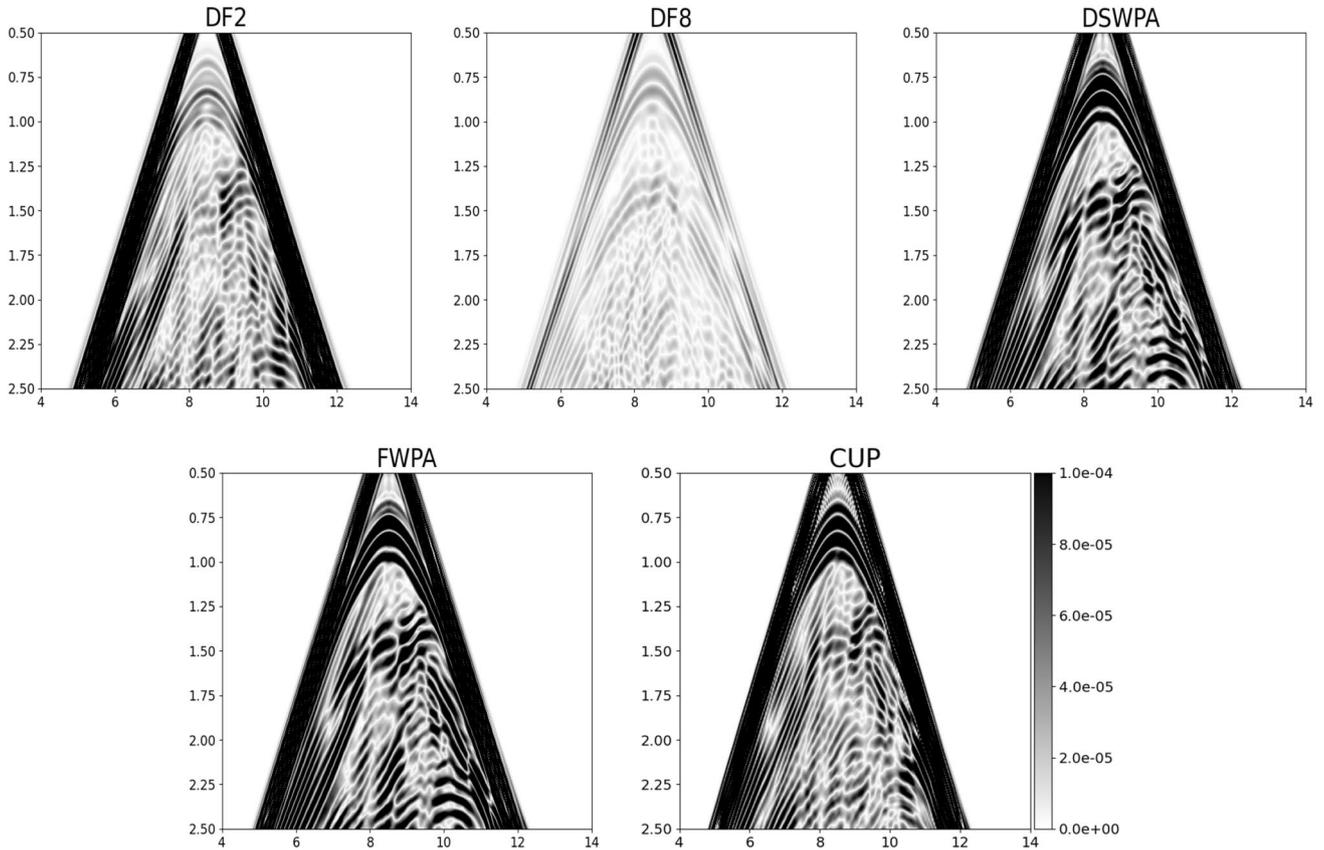

**Fig. 16** **Test case 3 (Marmousi).** Absolute difference between the reference seismogram and the seismogram obtained through the implementation of the DF2, DF8, DSWPA, FWPA, and CUP methods, with $\Delta x = \Delta y = 20$ m at $t = 2.5$ s

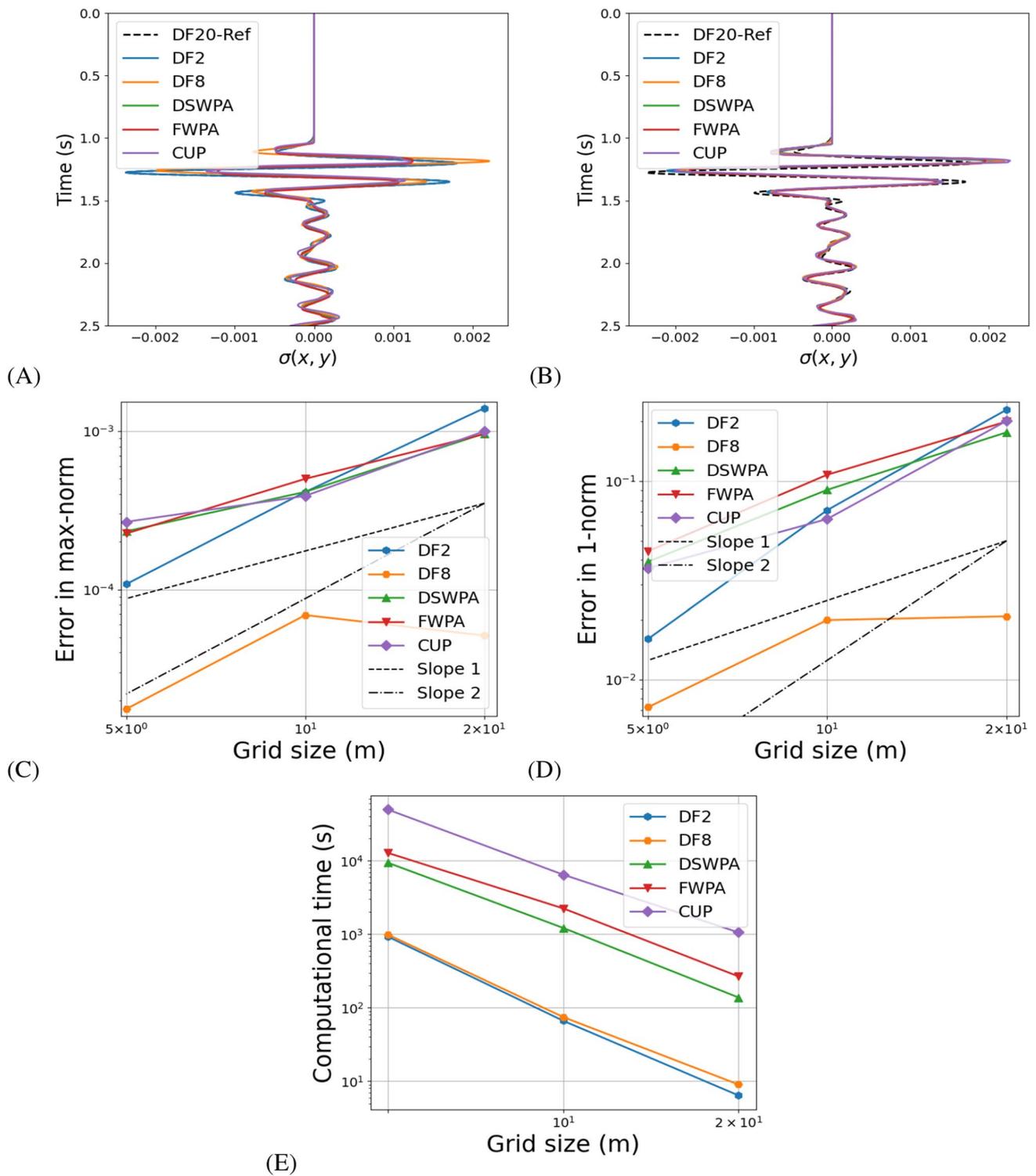

**Fig. 17 Test case 3 (Marmousi).** Top row: a comparison with the reference seismogram recorded at the receiver located at $x = 7000$ m and $t = 2.5$ s for the DF2, DF8, DSWPA, FWPA, and CUP methods with $\Delta x = \Delta y = 20$ m (**A**) and $\Delta x = \Delta y = 10$ m (**B**). Bottom row: log-log plots show the error versus grid size for the same methods using the max-norm (**C**), 1-norm (**D**), and a plot of grid size versus computational time for this test case (**E**). The lines labeled slope 1 and slope 2 represent references for 1st and 2nd order of convergence, respectively

**Fig. 18** **Test case 4.** A typical velocity field of Santos Basin. The domain spans from $x_0 = 20{,}000$ m to $x_1 = 50{,}000$ m and $y_0 = 0$ to $y_1 = 8000$ m

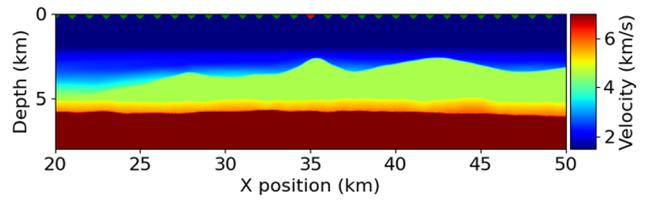

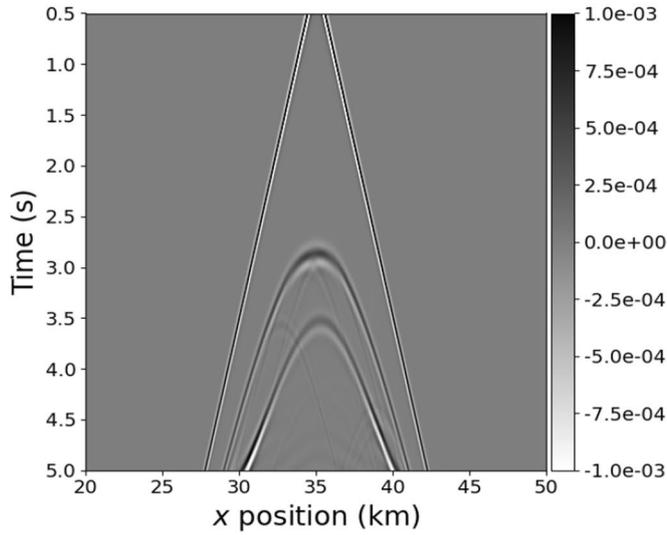

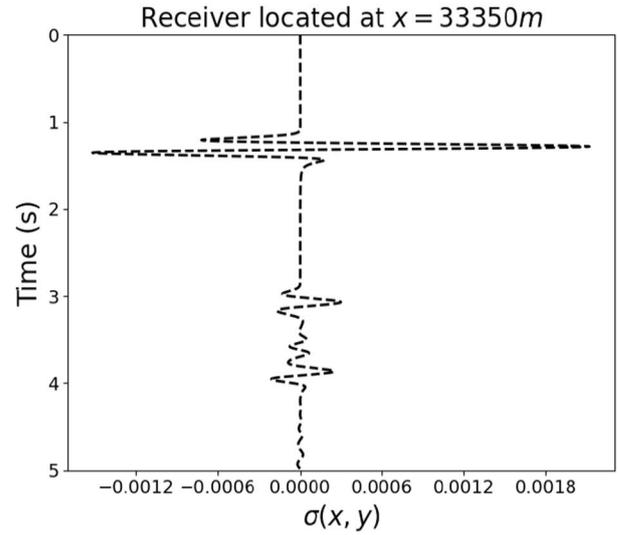

**Fig. 19** **Test case 4.** Reference seismogram of the $\sigma(x, y)$ field, computed using a DF20 method on a mesh consisting of 9601 points in $x$ and 4001 in $y$ directions. The spatial resolution was $\Delta x = 3.125$ m and $\Delta y = 2$ m, and the simulation was conducted until a final time of $t = 5$ s, with CFL number set to 0.25

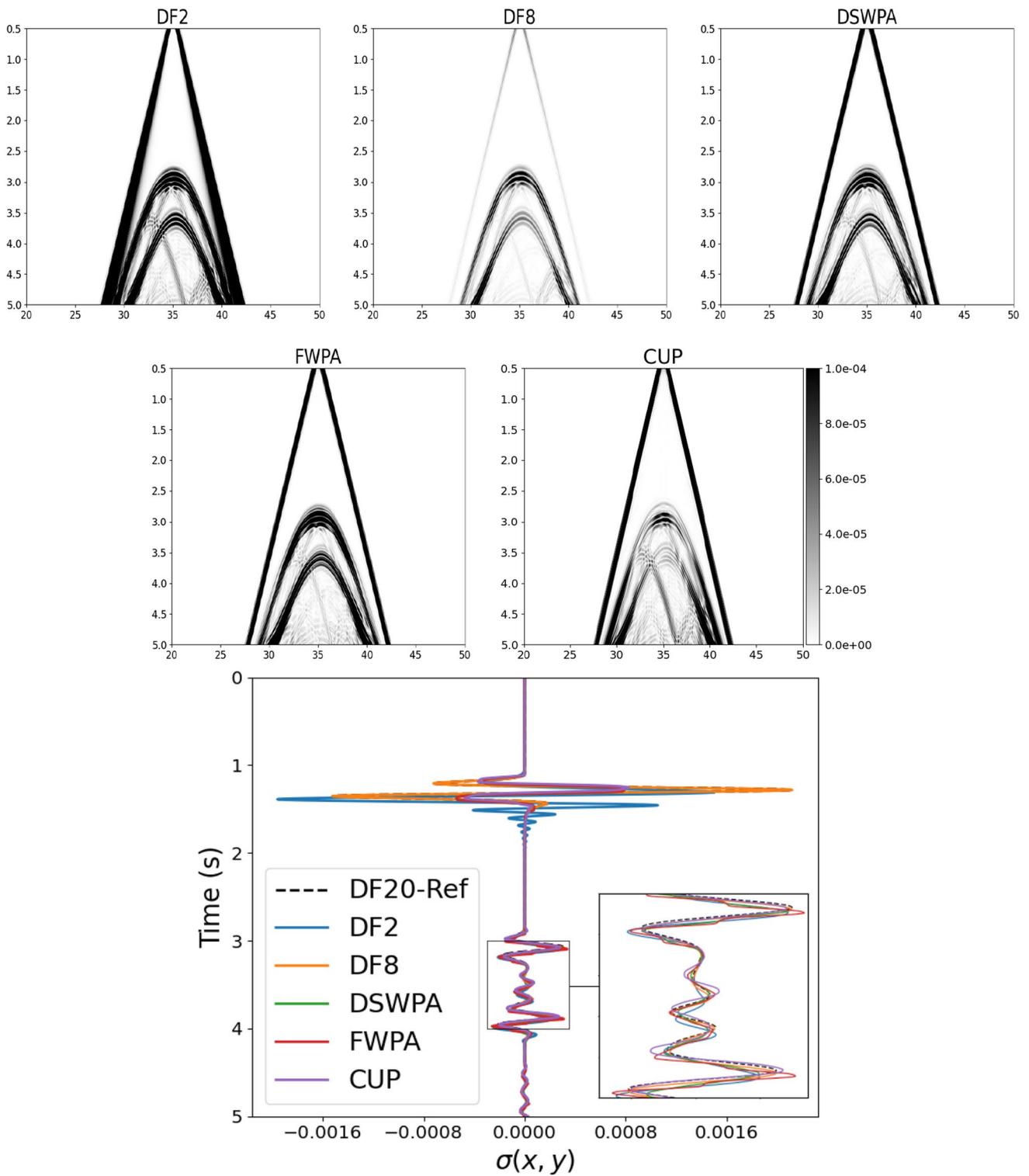

**Fig. 20 Test case 4.** The first five plots display the absolute difference between the reference seismogram and the seismograms obtained using the DF2, DF8, DSWPA, FWPA, and CUP methods, using a spatial resolution of $\Delta x = 25$ m and $\Delta y = 16$ m at time $t = 5$ s. The last plot illustrates the comparison with the reference seismogram recorded for the receiver located at $x = 33{,}350$ m

**Author Contributions** **Juan Barrios:** Writing—original draft, Visualization, Validation, Software, Methodology, Investigation, Formal analysis, Conceptualization, Writing—review & editing. **Pedro S. Peixoto:** Writing—review & editing, Supervision, Resources, Project administration, Methodology, Funding acquisition, Conceptualization, Writing—original draft. **Felipe A. G. Silva:** Conceptualization, Methodology, Resources, Software, Supervision, Writing—review & editing. The corresponding author has read the journal policies and submit this manuscript in accordance with those policies.

**Funding**   Juan Barrios acknowledge the financial support provided by the Agência Nacional do Petróleo, Gás Natural e Biocombustíveis under grant ANP20714-2 STMI-Shell and partial support from Coordenação de Aperfeiçoamento de Pessoal de Nível Superior - Brasil (CAPES) - Finance Code 001. Pedro S. Peixoto acknowledge the financial support provided by Agência Nacional do Petróleo, Gás Natural e Biocombustíveis under grants ANP20714-2 STMI-Shell, by São Paulo Research Foundation (FAPESP) under grant 2021/06176-0 and by the National Council for Scientific and Technological Development (CNPq) under grant 303436/2022-0. Felipe A. G. Silva acknowledge the financial support provided by the Agência Nacional do Petróleo, Gás Natural e Biocombustíveis under grant ANP20714-2 STMI-Shell.

**Data Availability**   The Santos Basin model data can be shared upon request from the corresponding author. The SEG/EAGE and Marmousi data used and all other data generated and analyzed during this study can be obtained using the codes openly provided in the GitHub repository (https://github.com/jbarrios94/ElasticWavePropagation.git). All figures and results are reproducible from the codes in the same repository.

## Declarations

**Ethics approval and consent to participate**   Not applicable.

**Consent for publication**   Not applicable.

**Competing interests**   The authors declare that they have no competing interests as defined by Discover, or other interests that might be perceived to influence the results and/or discussion reported in this paper.